\documentclass[dvips,11pt]{article}
\usepackage[latin1]{inputenc}
\usepackage{amsmath,amsthm,amsfonts}
\usepackage[numbers]{natbib}
\usepackage{enumerate}
\usepackage[ps2pdf,colorlinks]{hyperref}

\theoremstyle{plain}
\newtheorem{theorem}{Theorem}[section]
\newtheorem{lemma}[theorem]{Lemma}

\theoremstyle{definition}
\newtheorem{example}[theorem]{Example}
\newtheorem{remark}[theorem]{Remark}

\newcommand{\reff}[1]{(\ref{#1})}
\newcommand{\mailto}[1]{\href{mailto:#1}{\nolinkurl{#1}}}

\newcommand{\Z}{\mathbb{Z}}

\newcommand{\R}{\mathbb{R}}

\newcommand{\mB}{\mathcal{B}}

\newcommand{\mS}{\mathcal{S}}

\newcommand{\E}{\operatorname{E}}
\newcommand{\pr}{\operatorname{P}}
\newcommand{\wto}{\stackrel{w}{\to}}

\newcommand{\XLB}{X^{\rm{LB}}}

\newcommand{\alphaLB}{\alpha^{\rm{LB}}}
\newcommand{\QLB}{Q^{\rm{LB}}}

\newcommand{\la}{\leftarrow}
\newcommand{\ra}{\rightarrow}
\newcommand{\prm}{\mathcal{P}}
\newcommand{\Sprod}{S_1 \times S_2}
\newcommand{\Sprodp}{S_1' \times S_2'}

\newcommand{\proj}{\pi}

\newcommand{\idmap}{\chi}

\newcommand{\rel}[1]{\stackrel{#1}{\sim}}
\newcommand{\stochrel}[1]{\stackrel{#1}{\sim_{\rm st}}}

\newcommand{\Rst}{R_{\rm{st}}}
\newcommand{\simst}{\sim_{\rm{st}}}
\newcommand{\eqst}{=_{\rm{st}}}
\newcommand{\lest}{\le_{\rm{st}}}

\newcommand{\pre}{\preceq}

\newcommand{\Rsum}{R^{\rm{sum}}}
\newcommand{\approxst}{\approx_{\rm{st}}}
\newcommand{\leM}{\le_M}
\newcommand{\maj}{\pre^{\rm m}}
\newcommand{\majst}{\pre^{\rm m}_{\rm st}}
\newcommand{\wmaj}{\pre^{\rm wm}}
\newcommand{\wmajst}{\pre^{\rm wm}_{\rm st}}

\begin{document}



\title{Stochastic relations of \\ random variables and processes}
\author{
 Lasse Leskelä\thanks{
 Postal address: Helsinki University of Technology,
 PO Box 1100, 02015 TKK, Finland.
 URL: \url{www.iki.fi/lsl/}.
 Email: \protect\mailto{lasse.leskela@iki.fi}}
}
\date{\today}
\maketitle

\begin{abstract}
This paper generalizes the notion of stochastic order to a relation between probability measures
over arbitrary measurable spaces. This generalization is motivated by the observation that for the
stochastic ordering of two stationary Markov processes, it suffices that the generators of the
processes preserve some, not necessarily reflexive or transitive, subrelation of the order
relation. The main contributions of the paper are: a functional characterization of stochastic
relations, necessary and sufficient conditions for the preservation of stochastic relations, and
an algorithm for finding subrelations preserved by probability kernels. The theory is illustrated
with applications to hidden Markov processes, population processes, and queueing systems.
\end{abstract}

%


\vspace{1ex}



\section{Introduction}
\label{sec:introduction}

Comparison techniques based on stochastic orders~\cite{muller2002,shaked2007,szekli1995} are key
to obtaining upper and lower bounds for complicated random variables and processes in terms of
simpler random elements. Consider for example two ergodic discrete-time Markov processes $X$ and
$Y$ with stationary distributions $\mu_X$ and $\mu_Y$, taking values in a common ordered state
space, and denote by $\lest$ the corresponding stochastic order. Then the upper bound
\begin{equation}
 \label{eq:stationaryOrder}
 \mu_X \lest \mu_Y
\end{equation}
can be established~\cite{kamae1977} without explicit knowledge of $\mu_X$ by verifying that the
corresponding transition probability kernels $P_X$ and $P_Y$ satisfy
\begin{equation}
 \label{eq:kernelOrder}
 x \le y \implies P_X(x,\cdot) \lest P_Y(y,\cdot).
\end{equation}
Analogous conditions for continuous-time Markov processes on countable spaces have been derived by
Whitt~\cite{whitt1986} and Massey~\cite{massey1987}, and later extended to more general jump
processes by Brandt and Last~\cite{brandt1994}.

The starting point of this paper is to generalize the notion of stochastic order by denoting $X
\simst Y$, if there exists a coupling $(\hat X, \hat Y)$ of $X$ and $Y$ such that $\hat X \sim
\hat Y$ almost surely, where $\sim$ denotes some relation between the state spaces of $X$ and $Y$.
The main motivation for this definition is that~\reff{eq:kernelOrder} is by no means necessary
for~\reff{eq:stationaryOrder}; a less stringent sufficient condition is that
\begin{equation}
 \label{eq:kernelRelation}
 x \sim y \implies P_X(x,\cdot) \simst P_Y(y,\cdot)
\end{equation}
for some, not necessarily symmetric or transitive, nontrivial subrelation of the underlying order
relation. Another advantage of the generalized definition is that $X$ and $Y$ are no longer
required to take values in the same state space, leading to greater flexibility in the search for
bounding random elements $Y$. For example, to study whether $f(X) \lest g(Y)$ for some given real
functions $f$ and $g$ defined on the state spaces of $X$ and $Y$, we may define a relation $x \sim
y$ by the condition $f(x) \le g(y)$ \cite{doisy2000}.

The main contributions of the paper are: a functional characterization of stochastic relations,
necessary and sufficient conditions for the preservation of stochastic relations in the sense
of~\reff{eq:kernelRelation}, and an algorithm for finding subrelations preserved by probability
kernels. The functional characterization (Section~\ref{sec:stochasticRelations}) is given in terms
of relational conjugates that were implicitly defined by Strassen~\cite[Theorem 11]{strassen1965},
and the proof goes along similar lines, the new feature being the use of compact sets and upper
semicontinuous functions instead of completions of measures. López and Sanz have characterized the
preservation of stochastic relations for Markov processes on countable spaces in terms of a subtle
order construction~\cite{lopez2002}. Section~\ref{sec:preservation} describes an equivalent,
considerably simpler characterization based on relational conjugates, together with an iterative
algorithm for finding the maximal subrelation of a given relation preserved by a pair of
probability kernels. The main results are extended to the context of general random processes and
Markov processes in Section~\ref{sec:processes}. Applications to hidden Markov processes,
population processes, and queueing systems are discussed in Section~\ref{sec:applications}.
Section~\ref{sec:conclusion} concludes the paper.

\section{Stochastic relations}
\label{sec:stochasticRelations}

\subsection{Definitions}
\label{sec:definitions}

Let $(S_1,\mS_1)$ and $(S_2,\mS_2)$ be measurable spaces, and denote by $\prm(S_i)$ the family of
probability measures on $(S_i,\mS_i)$. Unless otherwise mentioned, all spaces shall implicitly be
assumed Polish (complete separable metrizable) and equipped with the Borel sigma-algebra. A
\emph{coupling} of probability measures $\mu_1 \in \prm(S_1)$ and $\mu_2 \in \prm(S_2)$ is a
probability measure $\mu \in \prm(\Sprod)$ with marginals $\mu_1$ and $\mu_2$, that is, $\mu \circ
\proj_i^{-1} = \mu_i$ for $i=1,2$, where $\proj_i$ denotes the projection map from $\Sprod$ onto
$S_i$. If $\mu$ is a coupling of $\mu_1$ and $\mu_2$, we also say that $\mu$ couples $\mu_1$ and
$\mu_2$ \cite{lindvall1992,thorisson2000}.

A \emph{measurable relation} between $S_1$ and $S_2$ is measurable subset of $\Sprod$. All
relations in this paper are assumed to be closed (in the product topology of $\Sprod$), if not
otherwise mentioned. Given a nontrivial ($R \neq \emptyset$) measurable relation $R$ between $S_1$
and $S_2$, we write $x_1 \sim x_2$, if $(x_1,x_2) \in R$. For probability measures $\mu_1 \in
\prm(S_1)$ and $\mu_2 \in \prm(S_2)$ we denote
\[
\mu_1 \simst \mu_2,
\]
and say that $\mu_1$ is stochastically related to $\mu_2$, if there exists a coupling $\mu$ of
$\mu_1$ and $\mu_2$ such that $\mu(R) = 1$. The relation $\Rst = \{ (\mu_1,\mu_2) : \mu_1 \simst
\mu_2 \}$ is called the \emph{stochastic relation} generated by $R$. Observe that two Dirac
measures satisfy $\delta_{x_1} \simst \delta_{x_2}$ if and only if $x_1 \sim x_2$. In this way the
stochastic relation $\Rst$ may be regarded as a natural randomization of the underlying relation
$R$.

A random variable $X_1$ is \emph{stochastically related} to a random variable $X_2$, denoted by
$X_1 \simst X_2$, if the distribution of $X_1$ is stochastically related to the distribution of
$X_2$. Observe that $X_1$ and $X_2$ do not need to be defined on the same probability space.
Recall that a coupling of random variables $X_1$ and $X_2$ is a bivariate random variable whose
distribution couples the distributions of $X_1$ and $X_2$. Hence $X_1 \simst X_2$ if and only if
there exists a coupling $(\hat X_1, \hat X_2)$ of $X_1$ and $X_2$ such that $\hat X_1 \sim \hat
X_2$ almost surely.

\begin{example}[Stochastic equality]
\label{exa:equality}
The stochastic relation generated by the equality relation $\{(x,y): x=y\}$ on $S$ is the equality
on $\prm(S)$. Hence $X \eqst Y$ if and only if $X$ and $Y$ have the same distribution.
\end{example}

\begin{example}[Stochastic $\epsilon$-distance]
Define a relation on the real line by denoting $x \approx y$, if $|x-y| \le \epsilon$. If $X_1
\approxst X_2$, then the cumulative distribution functions of $X_1$ and $X_2$ satisfy
\begin{equation}
 \label{eq:realLine}
 F_2(x-\epsilon) \le F_1(x) \le F_2(x + \epsilon) \quad \text{for all $x$}.
\end{equation}
Conversely, if~\reff{eq:realLine} holds, it is not hard to verify that the quantile functions
$G_i(r) = \inf \{x: F_i(x) \ge r\}$ satisfy $| G_1(r) - G_2(r)| \le \epsilon$ for all $r \in
(0,1)$. Hence the bivariate random variable $\hat X = (G_1(\xi), G_2(\xi))$, with $\xi$ uniformly
distributed on $(0,1)$, couples $X_1$ and $X_2$ and satisfies $\hat X_1 \approx \hat X_2$ with
probability one. Thus~\reff{eq:realLine} is necessary and sufficient for $X_1 \approxst X_2$.
\end{example}

\begin{example}[Stochastic majorization]
\label{exa:majorization} Let $S$ be closed subset of $\R^n$, and denote by $x_{[1]} \ge \cdots \ge
x_{[n]}$ the components of $x \in S$ in decreasing order. The \emph{weak majorization order} on
$S$ is defined by denoting $x \wmaj y$, if $\sum_{i=1}^k x_{[i]} \le \sum_{i=1}^k y_{[i]}$ for all
$k=1,\dots,n$; and the \emph{majorization order} by denoting $x \maj y$, if $x \wmaj y$ and
$\sum_{i=1}^n x_i = \sum_{i=1}^n y_i$. The $\maj$-increasing real functions are called
Schur-convex, and a function is $\wmaj$-increasing if and only if it is coordinatewise increasing
and Schur-convex~\cite[Theorem 3.A.8]{marshall1979}. The standard characterization for stochastic
orders (Remark~\ref{rem:stochasticOrder} in Section~\ref{sec:characterization}) hence shows that
$X \majst Y$ (resp.\ $X \wmajst Y$) if and only if $\E f(X) \le \E f(Y)$ for all positive
measurable Schur-convex (resp. coordinatewise increasing Schur-convex) functions $f$.
\end{example}

\subsection{Relational conjugates}
\label{sec:conjugates}

To develop a convenient way to check whether two probability measures are stochastically related
or not, we shall define the \emph{right conjugate} of $B_1 \subset S_1$ and the \emph{left
conjugate} of $B_2 \subset S_2$ with respect to a relation $R$ by
\begin{align*}
 B_1^\ra &= \cup_{x \in B_1} \{ y \in S_2: x\sim y\}, \\
 B_2^\la &= \cup_{y \in B_2} \{ x \in S_1: x\sim y\}.
\end{align*}
The conjugates of positive functions $f_i$ on $S_i$ are defined analogously by
\begin{align*}
 f_1^\ra(y) &= \sup_{x\in S_1: x\sim y} f_1(x), \quad y \in S_2, \\
 f_2^\la(x) &= \sup_{y\in S_2: x\sim y} f_2(y), \quad x \in S_1,
\end{align*}
where we adopt the convention that the supremum of the empty set is zero. Relational conjugates of
sets and functions are interlinked via
\begin{equation}
\label{eq:dualitySets}
 (1_{B_1})^\ra = 1_{B_1^\ra},
\end{equation}
where $1_{B_1}$ denotes the indicator function of $B_1$, and
\begin{equation}
\label{eq:dualityFunctions}
\{x: f_1(x) > r\}^\ra = \{y: f_1^\ra(y) > r\},
\end{equation}
which is valid for all $r \ge 0$.

The following result summarizes the basic topological properties of right conjugates. By symmetry,
analogous results are valid for left conjugates.

\begin{lemma}
\label{the:conjugatesTopological}
Let $R$ be a closed relation between two Polish spaces. Then:
\begin{enumerate}[(i)]
\item \label{pro:conjugateCompact}
$B^\ra$ is closed for compact $B$. Especially, $\{x\}^\ra$ is closed for all $x$.
\item \label{pro:conjugateUSC}
$f^\ra$ is upper semicontinuous (u.s.c.) for all positive u.s.c.\ $f$ on $S_1$ with compact
support.
\end{enumerate}
\end{lemma}
\begin{proof}
Assume $B$ is compact, and consider a sequence $y_n \to y$ such that $y_n \in B^\ra$ for all $n$.
Then for all $n$ there exists $x_n \in B$ such that $x_n \sim y_n$. Because $B$ is compact, there
exists $x \in B$ such that $x_n \to x$ as $n \to \infty$ along some subsequence of the natural
numbers. Hence $(x_n,y_n) \to (x,y)$ as $n \to \infty$ along the same subsequence, which implies
that $x \sim y$, and thus $y \in B^\ra$.

Assume next that $f$ is positive u.s.c.\ with compact support on $S_1$. We shall first show that
\begin{equation}
 \label{eq:dualityClosed}
 \{x: f(x) \ge r \}^\ra = \{y: f^\ra(y) \ge r\} \quad \text{for all} \ r > 0.
\end{equation}
Observe first that if $y$ belongs to the left side of~\reff{eq:dualityClosed}, then $f(x) \ge r$
for some $x \sim y$, so that $f^\ra(y) \ge r$. To prove the converse statement, assume next that
$f^\ra(y) \ge r$. Then the sets $K_n = \{f \ge r - 1/n\} \cap \{y\}^\la$ are nonempty and compact
for all $n > 1/r$, because $\{y\}^\la$ is closed by property~\reff{pro:conjugateCompact}. Hence
Cantor's intersection theorem implies that
\[
 \{x: f(x) \ge r\} \cap \{y\}^\la = \cap_{n>1/r} K_n,
\]
is nonempty, so that $f(x) \ge r$ for some $x \sim y$. We may now use~\reff{eq:dualityClosed}
together with property~\reff{pro:conjugateCompact} to conclude that $\{y: f^\ra(y) \ge r\}$ is
closed for all $r>0$. Obviously, $\{y: f^\ra(y) \ge 0\} = S_2$ is closed as well.
\end{proof}

\subsection{Functional characterization}
\label{sec:characterization}

The following result characterizes sto\-chas\-tic relations using relational conjugates of sets
and functions. The key part of the characterization is essentially Strassen's
Theorem~11~\cite{strassen1965}, written in a new notation. The new contributions
are~\reff{pro:testSetsCompact} and~\reff{pro:testFunctionsUSC}, providing classes of test sets and
functions with Borel-measurable conjugates (Lemma~\ref{the:conjugatesTopological}) that are large
enough to characterize stochastic relations without resorting to completions of measures.

\begin{theorem}
\label{the:relationTest}
Let $R$ be a closed relation between Polish spaces $S_1$ and $S_2$. Then $\mu \simst \nu$ is
equivalent to each of the following:
\begin{enumerate}[(i)]
\item \label{pro:testSetsMeasurable} $\mu(B) \le \nu(B^\ra)$ for all measurable $B$ such that $B^\ra$ is measurable.
\item \label{pro:testSetsCompact} $\mu(B) \le \nu(B^\ra)$ for all compact $B$.
\item \label{pro:testFunctionsMeasurable} $\int_{S_1} f \, d\mu \le \int_{S_2} f^\ra \, d\nu$ for all positive measurable $f$
such that $f^\ra$ is measurable.
\item \label{pro:testFunctionsUSC} $\int_{S_1} f \, d\mu \le \int_{S_2} f^\ra \, d\nu$ for all
positive u.s.c.\ $f$ with compact support.
\end{enumerate}
\end{theorem}

\begin{remark}
\label{rem:stochasticOrder} If $R$ is an order (reflexive and transitive) relation on $S$, then
using the properties $B \subset B^\ra = (B^\ra)^\ra$ and $f \le f^\ra = (f^\ra)^\ra$ we see that
\reff{pro:testSetsMeasurable} and \reff{pro:testFunctionsMeasurable} in
Theorem~\ref{the:relationTest} become equivalent to well-known characterizations of stochastic
orders~\cite{kamae1977,strassen1965}:
\begin{itemize}
\item[(i')] $\mu(B) \le \nu(B)$ for all measurable upper sets $B$.
\item[(iii')] $\int_S f \, d\mu \le \int_S f \, d\nu$ for all measurable positive increasing functions $f$.
\end{itemize}
\end{remark}

\begin{remark}
When $S_1$ and $S_2$ are countable, the measurability requirements of
Theorem~\ref{the:relationTest} become void, and the word "compact" becomes replaced by "finite".
\end{remark}

\begin{proof}[Proof of Theorem~\ref{the:relationTest}]
$\mu \simst \nu$ $\implies$ \reff{pro:testSetsMeasurable}. Let $\lambda$ be a coupling of $\mu$
and $\nu$ such that $\lambda(R) = 1$. Then because
\[
(B \times S_2) \cap R = (B \times B^\ra) \cap R,
\]
we see that
\[
\mu(B) = \lambda(B \times S_2) = \lambda(B \times B^\ra) \le \lambda(S_1 \times B^\ra) =
\nu(B^\ra)
\]
for all measurable $B \subset S_1$ such that $B^\ra$ is measurable.

\reff{pro:testSetsMeasurable} $\implies$ \reff{pro:testSetsCompact}. Clear by
Lemma~\ref{the:conjugatesTopological}.

\reff{pro:testSetsCompact} $\implies$ \reff{pro:testFunctionsUSC}. Let $f$ be a positive compactly
supported u.s.c.\ function on $S_1$. Then equality~\reff{eq:dualityClosed} shows that
\[
 \mu(\{x: f(x) \ge r\}) \le \nu(\{y: f^\ra(y) \ge r\})
\]
for all $r > 0$. The validity of~\reff{pro:testFunctionsUSC} hence follows by integrating both
sides of the above inequality with respect to $r$ over $(0,\infty)$.

\reff{pro:testFunctionsUSC} $\implies$ $\mu \simst \nu$. By virtue of~\cite[Theorem
7]{strassen1965}, it suffices to show that
\begin{equation}
\label{eq:strassenNew}
 \int_{S_1} f \, d\mu + \int_{S_2} g \, d\nu \le \sup_{(x,y) \in R} (f(x) + g(y))
\end{equation}
for all bounded continuous $f$ and $g$ on $S_1$ and $S_2$, respectively, and without loss of
generality we may assume $f$ and $g$ are positive and bounded by one. Given any such functions $f$
and $g$, and a number $\epsilon > 0$, choose a compact set $K \subset S_1$ such that $\mu(K^c) \le
\epsilon$ \cite[Theorem 1.3]{billingsley1999}, and define $f_0 = f 1_K$. Because $f_0$ is u.s.c.,
we see using~\reff{pro:testFunctionsUSC} that $\int_{S_1} f_0 \, d\mu \le \int_{S_2} f_0^\ra \,
d\nu$, so that
\begin{equation}
 \label{eq:strassenStep}
 \int_{S_1} f \, d\mu + \int_{S_2} g \, d\nu
 \le \int_{S_2} (f_0^\ra + g) \, d\nu + \epsilon.
\end{equation}
In light of~\reff{eq:dualitySets}, assumption~\reff{pro:testFunctionsUSC} further implies that
$\mu(K) \le \nu(K^\ra)$, because $1_K$ is u.s.c.. Thus $\nu((K^\ra)^c) \le \epsilon$, so by
splitting the $\nu$-integral into $K^\ra$ and its complement we see that
\[
 \int (f_0^\ra + g) \, d\nu \le \sup_{y \in K^\ra} (f_0^\ra(y) + g(y)) + 2 \epsilon.
\]
Because $f_0^\ra \le f^\ra$ and $K^\ra \subset S_1^\ra$, the above inequality combined
with~\reff{eq:strassenStep} shows that
\[
\int_{S_1} f \, d\mu + \int_{S_2} g \, d\nu \le \sup_{y \in S_1^\ra} (f^\ra(y) + g(y)) + 3
\epsilon.
\]
After letting $\epsilon \to 0$ and observing that
\[
\sup_{y \in S_1^\ra} (f^\ra(y) + g(y)) = \sup_{(x,y) \in R} (f(x) + g(y)),
\]
we may conclude that~\reff{eq:strassenNew} holds.

Finally, observe that the proof of \reff{pro:testSetsMeasurable} $\implies$
\reff{pro:testFunctionsMeasurable} is completely analogous to the proof of
\reff{pro:testSetsCompact} $\implies$ \reff{pro:testFunctionsUSC}, and the implication
\reff{pro:testFunctionsMeasurable} $\implies$ \reff{pro:testFunctionsUSC} follows immediately by
Lemma~\ref{the:conjugatesTopological}.
\end{proof}

\section{Preservation of stochastic relations}
\label{sec:preservation}

\subsection{Coupling of probability kernels}

Monotone functions are key objects in the study of order relations. When passing from orders to
general relations, the role of monotone functions is taken over by function pairs $(f_1,f_2)$ such
that $x_1 \sim x_2 \implies f_1(x_1) \sim f_2(x_2)$. To study stochastic relations, we need a
randomized version of the above property. Recall that a \emph{probability kernel} from a
measurable space $S$ to a measurable space $S'$ is a mapping $P: S \times \mS' \to \R$ such that
$P(x,\cdot)$ is a probability measure for all $x$, and $x \mapsto P(x,B)$ is measurable for all $B
\in \mS'$. Probability kernels may alternatively be viewed as mappings $\prm(S) \ni \mu \mapsto
\mu P \in \prm(S')$ by defining $\mu P(B) = \int_S P(x,B) \, \mu(dx)$.

Given a closed relation $R$ between Polish spaces $S_1$ and $S_2$, and probability kernels $P_1$
on $S_1$ and $P_2$ on $S_2$, we say that the pair $(P_1,P_2)$ \emph{stochastically preserves} $R$,
if any of the equivalent conditions in Theorem~\ref{the:preservation} holds.

\begin{theorem}
\label{the:preservation} The following are equivalent:
\begin{enumerate}[(i)]
\item \label{pro:preservation1} $x_1 \sim x_2 \implies P_1(x_1,\cdot) \simst P_2(x_2,\cdot)$.
\item \label{pro:preservation2} $\mu_1 \simst \mu_2 \implies \mu_1 P_1 \simst \mu_2 P_2$.
\item \label{pro:preservation3} $P_1(x_1,B) \le P_2(x_2,B^\ra)$ for all $x_1 \sim x_2$ and compact $B \subset S_1$.
\end{enumerate}
\end{theorem}

\begin{proof}
\reff{pro:preservation1} $\implies$ \reff{pro:preservation2}. Given $\mu_1 \simst \mu_2$, choose a
coupling of $\mu_1$ and $\mu_2$ such that $\mu(R)=1$. Theorem~\ref{the:relationTest} then shows
that
\[
 \mu_1 P_1(B)
 = \int_R P_1(x_1,B) \, \mu(dx)
 \le \int_R P_2(x_2,B^\ra) \, \mu(dx)
 = \mu_2 P_2(B^\ra)
\]
for all compact $B \subset S_1$, so that $\mu_1 P_1 \simst \mu_2 P_2$.

The implication \reff{pro:preservation2} $\implies$ \reff{pro:preservation1} follows immediately
by choosing $\mu_i = \delta_{x_i}$, while the equivalence \reff{pro:preservation1}
$\Longleftrightarrow$ \reff{pro:preservation3} is clear by Theorem~\ref{the:relationTest}.
\end{proof}

\begin{remark}
A probability kernel $P$ is said to stochastically preserve a relation $R$ on $S$, if $x_1 \sim
x_2 \implies P(x_1,\cdot) \simst P(x_2,\cdot)$. Order-preserving probability kernels are usually
called monotone~\cite{muller2002}.
\end{remark}

The main result of this section is the following coupling characterization of relation-preserving
pairs of probability kernels. For technical reasons related to local uniformization of Markov jump
processes in Section~\ref{sec:markov}, we shall consider probability kernels $P_i$ from $S_i$ to
$S_i'$, where $S_i'$ is a measurable space not necessarily equal to $S_i$. A probability kernel
$P$ from $\Sprod$ to $S_1' \times S_2'$ is called a \emph{coupling} of probability kernels $P_1$
and $P_2$, if the probability measure $P(x,\cdot)$ couples the probability measures
$P_1(x_1,\cdot)$ and $P_2(x_2,\cdot)$ for all $x = (x_1,x_2) \in \Sprod$.

\begin{theorem}
\label{the:kernelCouplingGeneral}
Given closed relations $R$ between $S_1$ and $S_2$, and $R'$ between $S_1'$ and $S_2'$, assume
that
\[
x_1 \rel{R} x_2 \implies P_1(x_1,\cdot) \stochrel{R'} P_2(x_2,\cdot).
\]
Then there exists a coupling $P$ of $P_1$ and $P_2$ such that $P(x,R') = 1$ for all $x \in R$.
\end{theorem}

The proof of Theorem~\ref{the:kernelCouplingGeneral} requires some preliminaries on topology and
measure theory that are discussed next. Denote by $\proj_i$ the projection from $S_1 \times S_2$
to $S_i$, and define the projection maps $\hat\proj_i: \prm(S_1 \times S_2) \to \prm(S_i)$ by
$\hat\proj_i \mu = \mu \circ \proj_i^{-1}$, so that $\hat\proj_i \mu$ equals the $i$-th marginal
of $\mu$. From now on, all sets of probability measures shall be considered as topological spaces
equipped with the weak topology.

\begin{lemma}
\label{the:continuousProjection}
The projection $\hat\proj_i: \prm(\Sprod) \to \prm(S_i)$ is continuous and open with respect to
the weak topology, $i=1,2$.
\end{lemma}
\begin{proof}
Assume that $\mu^n \wto \mu$ in $\prm(\Sprod)$, and let $f \in C_b(S_i)$. Then because $f \circ
\proj_i \in C_b(\Sprod)$, it follows that $\hat\proj_i \mu^n(f) = \mu^n(f \circ \proj_i) \to \mu(f
\circ \proj_i) = \hat\proj_i \mu(f)$. Hence $\hat \proj_i$ is continuous. The openness of
$\hat\proj_i$ follows from Eifler~\cite[Theorem 2.5]{eifler1976}, because the map $\proj_i$ is
continuous, open, and onto.
\end{proof}

\begin{lemma}
\label{the:compactCouplingSet}
For any $\mu_1 \in \prm(S_1)$ and $\mu_2 \in \prm(S_2)$, the set $K(\mu_1,\mu_2)$ of all couplings
of $\mu_1$ and $\mu_2$ is compact in the weak topology of $\prm(\Sprod)$.
\end{lemma}
\begin{proof}
Given $\epsilon > 0$, choose compacts sets $C_i \subset S_i$ such that $\mu_i(C_i^c) \le
\epsilon/2$, $i=1,2$. Define $C = C_1 \times C_2$. Then because $C^c = (C_1^c \times S_2) \cup
(S_1 \times C_2^c)$, it follows that
\[
 \mu(C^c) \le \mu_1(C_1^c) + \mu_2(C_2^c) \le \epsilon
\]
for all $\mu \in K(\mu_1,\mu_2)$. Hence $K(\mu_1,\mu_2)$ is relatively compact by Prohorov's
theorem. The equality $K(\mu_1,\mu_2) = \hat\proj_1^{-1}(\mu_1) \cap \hat\proj_2^{-1}(\mu_2)$
further shows that $K(\mu_1,\mu_2)$ is closed, because $\hat\proj_i$ are continuous by
Lemma~\ref{the:continuousProjection}.
\end{proof}

\begin{lemma}
\label{the:measurableKernel}
Let $P$ be a probability kernel from $S$ to $S'$. Then the map $x \mapsto P(x,\cdot)$ is
$\mB(S)/\mB(\prm(S'))$-measurable, where $\mB(\prm(S'))$ denotes the Borel $\sigma$-algebra
generated by the weak topology on $\prm(S')$.
\end{lemma}
\begin{proof}
For any $f \in C_b(S')$, one may check by approximating $f$ with simple functions that the map $x
\mapsto P(x,f)$ is measurable. Hence it follows that the set $\{x: P(x,\cdot) \in A\}$ is
$\mB(S)$-measurable for any $A = \cap_{k=1}^n \{ \mu: \mu(f_k) \in B_k \}$, where $f_k \in
C_b(S')$ and $B_k$ are open subsets of the real line. Because the sets of the above type form a
basis for the weak topology of $\prm(S')$, and because the space $\prm(S')$ equipped with the weak
topology is Polish~\cite{billingsley1999} and hence Lindelöf, it follows that any open set in
$\prm(S')$ can be represented as a countable union of the basis sets. Hence $\{x: P(x,\cdot) \in
A\}$ is measurable for all open subsets $A$ in $\prm(S')$, and the claim follows.
\end{proof}

A \emph{set-valued mapping} from a set $S$ to a set $S'$ is a function that assigns to each
element in $S$ a subset of $S'$. A set-valued mapping $F$ from a measurable space $S$ to a
topological space $S'$ is \emph{measurable}~\cite{wagner1977}, if the inverse image
\[
 F^-(A) = \{x \in S : F(x) \cap A \neq \emptyset \}
\]
is measurable for all closed $A \subset S'$.

\begin{lemma}
\label{the:setValuedMapping}
Let $P_i$ be probability kernels from $S_i$ to $S_i'$, $i=1,2$. Then the set-valued mapping  $F: x
\mapsto K(P(x_1,\cdot),P_2(x_2,\cdot))$ is measurable.
\end{lemma}
\begin{proof}
Because the set $F(x) \subset \prm(\Sprodp)$ is compact for all $x$ by
Lemma~\ref{the:compactCouplingSet}, it is sufficient to verify that $F^-(A)$ is measurable for all
open sets $A$ (Himmelberg~\cite[Theorem 3.1]{himmelberg1975}). Let us hence assume that $A \subset
\prm(\Sprodp)$ is open. Observe that $F^-(A) = \proj_1^{-1}(B_1) \cap \proj_2^{-1}(B_2)$, where
\[
B_i = \{x_i \in S_i: P_i(x_i,\cdot) \in \hat\proj_i(A)\}.
\]
Now Lemma~\ref{the:continuousProjection} implies that $\hat\proj_i(A)$ is open, and
Lemma~\ref{the:measurableKernel} further shows that $B_i$ is measurable. Thus $F^-(A)$ is
measurable.
\end{proof}

\begin{proof}[Proof of Theorem~\ref{the:kernelCouplingGeneral}]
Assume without loss of generality that $R' \neq \emptyset$, and let $\prm(R) = \{\mu \in
\prm(\Sprodp): \mu(R') = 1\}$. Because $R'$ is closed, it follows from Portmanteau's theorem that
$\prm(R')$ is closed. Define the set-valued mappings $F$ and $G$ from $\Sprod$ to $\prm(\Sprodp)$
by $F(x) = K(P_1(x_1,\cdot), P_2(x_2,\cdot))$, and
\[
 G(x) = \left\{
 \begin{aligned}
  F(x) \cap \prm(R'), \quad & x \in R, \\
  F(x),              \quad & \text{else}.
 \end{aligned}
\right.
\]
Then for any closed $A' \subset \prm(\Sprodp)$,
\[
G^-(A') = \left( R \cap F^-(\prm(R') \cap A') \right) \, \cup \, \left( R^c \cap F^-(A') \right).
\]
Because $\prm(R')$ is closed, it follows from Lemma~\ref{the:setValuedMapping} that $G^-(A')$ is
measurable. Moreover, because $F$ is compact-valued by Lemma~\ref{the:compactCouplingSet}, we may
conclude that $G$ is a measurable set-valued mapping such that $G(x)$ is compact and nonempty for
all $x$. A measurable selection theorem of Kuratowski and Ryll--Nardzewski~\cite{kuratowski1965}
(see alternatively Srivastava~\cite[Theorem 5.2.1]{srivastava1998}) now shows that there exists a
measurable function $g: \Sprod \to \prm(\Sprodp)$ such that $g(x) \in G(x)$ for all $x$. By
defining $P(x,B) = [g(x)](B)$ for $x \in \Sprod$ and measurable $B \subset \Sprodp$, we see that
$P$ is a probability kernel from $\Sprod$ to $\Sprodp$ with the desired properties.
\end{proof}

\subsection{Subrelation algorithm}
\label{sec:subrelation}

This section presents an algorithm for finding the maximal subrelation of a closed relation that
is stochastically preserved by a pair $(P_1,P_2)$ of continuous\footnote{A probability kernel $P$
from $S$ to $S'$ is called \emph{continuous} if $P(x_n, \cdot) \to P(x, \cdot)$ in distribution
whenever $x_n \to x$. In other words, $P$ is continuous if and only if the map $x \mapsto
P(x,\cdot)$ from $S$ to $\prm(S')$ is continuous, when $\prm(S')$ is equipped with the weak
topology.} probability kernels. Given a closed relation $R$ and continuous probability kernels
$P_i$ on $S_i$, $i=1,2$, define recursively the relations $R^{(n)}$ by $R^{(0)} = R$,
\[
 R^{(n+1)} = \left\{x \in R^{(n)}: \ (P_1(x_1,\cdot), P_2(x_2,\cdot)) \in \Rst^{(n)} \right\},
\]
and denote
\begin{equation}
 \label{eq:RStar}
 R^* = \bigcap_{n=0}^\infty R^{(n)}.
\end{equation}
Lemma~\ref{the:MDT} below shows that the relations $R^{(n)}$ are closed and hence measurable, so
the stochastic relations $\Rst^{(n)}$ are well-defined. The following theorem underlines the key
role of $R^*$ in characterizing the existence of subrelations stochastically preserved by a pair
of probability kernels.

\begin{theorem}
\label{the:subrelationKernel} Assume $P_1$ and $P_2$ are continuous. Then $R^*$ is the maximal
closed subrelation of $R$ that is stochastically preserved by $(P_1,P_2)$. Especially, there
exists a nontrivial closed subrelation stochastically preserved by $(P_1,P_2)$ if and only if $R^*
\neq \emptyset$.
\end{theorem}

The following three lemmas summarize the topological preliminaries required for the proof of
Theorem~\ref{the:subrelationKernel}.

\begin{lemma}
\label{the:closed} Let $R$ be a closed relation between Polish spaces $S_1$ and $S_2$. Then the
relation $\Rst$ is closed in the weak topology of $\prm(S_1) \times \prm(S_2)$.
\end{lemma}
\begin{proof} Assume $\mu_n \wto \mu$ and $\nu_n \wto \nu$ such that $\mu_n \simst \nu_n$ for all
$n$. Then for all $n$ there exists a coupling $\lambda_n$ of $\mu_n$ and $\nu_n$ such that
$\lambda_n(R) = 1$. Because the sequences $\mu_n$ and $\nu_n$ are tight, and because $\lambda_n (
(C_1 \times C_2)^c ) \le \mu_n(C_1^c) + \nu_n(C_2^c)$ for all compact $C_1$ and $C_2$, it follows
that the sequence $\lambda_n$ is tight, so there exists $\lambda \in \prm(\Sprod)$ such that
$\lambda_n \wto \lambda$ as $n \to \infty$ along some subsequence~\cite[Theorem
16.3]{kallenberg2002}. The continuity of $\hat \pi_i$ (Lemma~\ref{the:continuousProjection})
implies that $\lambda$ is a coupling of $\mu$ and $\nu$, and Portmanteau's theorem shows that
$\lambda(R) = 1$. Hence $\mu \simst \nu$.
\end{proof}

\begin{lemma}
\label{the:MDT} Given continuous probability kernels $P_1$ and $P_2$, define
\begin{equation}
\label{eq:defM} M(R) = \left\{x \in R: \ (P_1(x_1,\cdot), P_2(x_2,\cdot)) \in \Rst \right\}
\end{equation}
for measurable relations $R$. Then:
\begin{enumerate}[(i)]
\item \label{pro:monotoneM} $M(R) \subset M(R')$ for $R \subset R'$.
\item \label{pro:closedM}   $M$ maps closed relations into closed relations.
\end{enumerate}
\end{lemma}
\begin{proof}

For~\reff{pro:monotoneM} it suffices to observe that $R \subset R'$ implies $\Rst \subset \Rst'$.
For~\reff{pro:closedM}, observe that $M(R) = R \cap f^{-1}(\Rst)$, where the function $f: \Sprod
\to \prm(S_1) \times \prm(S_2)$ is defined by $f(x) = (P_1(x_1,\cdot), P_2(x_2,\cdot))$. Because
$\Rst$ is closed (Lemma~\ref{the:closed}) and $f$ is continuous, it follows that $M(R)$ is closed.
\end{proof}

\begin{lemma}
\label{the:relationConvergence}
Let $R^{(0)} \supset R^{(1)} \supset \cdots$ be closed relations between Polish spaces $S_1$ and
$S_2$, and let $R^* = \cap_{n=0}^\infty R^{(n)}$. Assume that $(\mu_1,\mu_2) \in \Rst^{(n)}$ for
all $n$. Then $(\mu_1, \mu_2) \in \Rst^*$.
\end{lemma}
\begin{proof}
By definition, for all $n$ there exists a coupling $\lambda_n$ of $\mu$ and $\nu$ such that
$\lambda_n(R^{(n)}) = 1$. Because the set of couplings of $\mu$ and $\nu$ is compact by
Lemma~\ref{the:compactCouplingSet}, there exists a coupling $\lambda$ of $\mu$ and $\nu$ such that
$\lambda_n \wto \lambda$ as $n \to \infty$ along a subsequence of $\Z_+$. Further, observe that
$\lambda_n(R^{(m)}) \ge \lambda_n(R^{(n)}) = 1$ for all $m \le n$, which implies that
$\lim_{n\to\infty} \lambda_n(R^{(m)}) = 1$ for all $m$. Portmanteau's theorem now shows that
$\lambda(R^{(m)}) = 1$ for all $m$, so it follows that $\lambda(R^*) = 1$.
\end{proof}

\begin{proof}[Proof of Theorem~\ref{the:subrelationKernel}]
Let $M$ be the map defined in~\reff{eq:defM}. If $x \in R^*$, then $x \in R^{(n+1)} = M(R^{(n)})$
shows that $(P_1(x_1, \cdot), P_2(x_2, \cdot)) \in \Rst^{(n)}$ for all $n$. Because the relations
$R^{(n)}$ are closed by Lemma~\ref{the:MDT}, we see using Lemma~\ref{the:relationConvergence} that
$(P_1(x_1, \cdot), P_2(x_2,\cdot)) \in \Rst^*$. Hence $(P_1,P_2)$ stochastically preserves $R^*$.
On the other hand, if $R'$ is a closed subrelation of $R$ that is stochastically preserved by
$(P_1,P_2)$, then $R' = M(R') \subset M(R) = R^{(1)}$ by Lemma~\ref{the:MDT}. Induction shows that
$R' \subset R^{(n)}$ for all $n$, and thus $R' \subset R^*$.
\end{proof}

\section{Random processes}
\label{sec:processes}

\subsection{Random sequences}
\label{sec:sequences} Given a relation $R$ between $S_1$ and $S_2$, the coordinatewise relation
between the product spaces $S_1^n$ and $S_2^n$ for $n \le \infty$ is defined by
\[
 R^n = \{(x,y) \in S_1^n \times S_2^n: x_i \sim y_i \ \text{for all} \ i\}.
\]
The stochastic relation generated by $R^n$ is called the \emph{stochastic coordinatewise
relation}. The following example shows that the stochastic coordinatewise relation of two random
sequences cannot be verified just by looking at the one-dimensional marginal distributions.
\begin{example}
Let $\xi_1$ and $\xi_2$ be independent random variables uniformly distributed on the unit
interval, and define $X = (\xi_1, \xi_1)$ and $Y = (\xi_1, \xi_2)$. Then $X_i \eqst Y_i$ for all
$i=1,2$, but $X$ and $Y$ are not related with respect to the stochastic coordinatewise equality.
\end{example}

The proof of the following result is a straightforward modification of its continuous-time
analogue Theorem~\ref{the:fidiCT}, and shall hence be omitted.
\begin{theorem}
\label{the:fidiDT}
Two random sequences $X$ and $Y$ satisfy $X \stochrel{R^n} Y$ if and only if $(X_{t_1}, \dots,
X_{t_k}) \stochrel{R^k} (Y_{t_1}, \dots, Y_{t_k})$ for all finite parameter combinations
$(t_1,\dots,t_k)$.
\end{theorem}

The following result, which is completely analogous to~\cite[Proposition~1]{kamae1977}, gives a
sufficient condition for $X \simst Y$ in terms of conditional probabilities. Let $P_{i}$ be a
probability kernel from $S_1^{i-1}$ to $S_1$ representing the regular conditional distribution of
$X_{i}$ given $(X_1,\dots,X_{i-1})$ , and define the kernels $Q_{i}$ in a similar way for $Y$
\cite[Theorem 6.3]{kallenberg2002}.

\begin{theorem}
\label{the:conditionalRelation}
Assume that $X_1 \simst Y_1$, and that for all $i$,
\[
P_{i}(x_1,\dots,x_{i-1}, dx_{i}) \simst Q_{i}(y_1,\dots,y_{i-1}, dy_{i})
\]
whenever $(x_1,\dots,x_{i-1}) \sim (y_1,\dots,y_{i-1})$. Then $X \simst Y$.
\end{theorem}
\begin{proof}
Let $\lambda_1$ be a coupling of the distributions of $X_1$ and $Y_1$ such that $\lambda_1(R) =
1$. For convenience, we shall use $x^n$ as a shorthand for $(x_1,\dots,x_n)$. By
Theorem~\ref{the:kernelCouplingGeneral} there exists for each $i$ a coupling $\Lambda_{i}$ of
probability kernels $P_{i}$ and $Q_{i}$ such that $\Lambda_i((x^{i-1},y^{i-1}),R) = 1$ whenever
$x^{i-1} \sim y^{i-1}$. Then it is easy to verify by induction that the probability measure
\[
 \lambda_n(B) = \idotsint 1(z^n \in B) \, \Lambda_n(z^{n-1},dz_n) \cdots \Lambda_2(z^1,dz_2) \, \lambda_1(dz_1)
\]
couples the distributions of $(X_1,\dots,X_n)$ and $(Y_1,\dots,Y_n)$, and $\lambda_n(R^n) = 1$ for
any finite $n$. In the case where $n$ is infinite, the proof is completed by applying
Theorem~\ref{the:fidiDT}.
\end{proof}

The next example shows that the condition in Theorem~\ref{the:conditionalRelation} is not
necessary in general.

\begin{example}
Let $X = (\xi, 1-\xi)$ and $Y = (2 \xi, 1-\xi)$, where $\xi$ has uniform distribution on the unit
interval. Then $X \lest Y$ by construction, but $P_2(x_1,\cdot) \lest Q_2(y_2,\cdot)$ only for
$x_1 \ge y_1/2$.
\end{example}

\subsection{Continuous-time random processes}
\label{sec:continuousTime}

Denote by $D_i = D_i(\R_+,S_i)$ the space of functions from $\R_+$ into $S_i$ that are
right-continuous and have left limits, and equip $D_i$ with the Skorohod topology, which makes it
Polish~\cite[Section 3.5]{ethier1986}. The coordinatewise relation between $D_1$ and $D_2$ is
defined by
\[
 R^D = \{(x,y) \in D_1 \times D_2: \ x(t) \sim y(t) \ \text{for all} \ t \in \R_+\},
\]
and we denote by $\Rst^D$ the corresponding stochastic relation between random processes with
paths in $D_i$ (identified as $D_i$-valued random elements). When there is no risk of confusion,
the same notation $\simst$ shall be used for a random process (corresponding to $\Rst^D$ and its
finite-dimensional distributions (corresponding to $\Rst^n$).

\begin{lemma}
\label{the:closedSkorohod}
$R^D$ is a closed relation between $D_1$ and $D_2$, whenever $R$ is closed.
\end{lemma}
\begin{proof}
Assume that $x_i$ and $x_i^n$ are functions in $D_i$ such that $x_i^n \to x_i$ as $n \to \infty$,
and $(x_1^n, x_2^n) \in R^D$ for all $n$. Denote $\Delta = \Delta_1 \cup \Delta_2$, where
$\Delta_i$ is the set of points $t \in \R_+$ where $x_i$ is discontinuous. It is well-known
\cite[Section 3.5]{ethier1986} that $\Delta_i$ is countable, and that $x_i^n(t) \to x_i(t)$ in
$S_i$ for all $t \notin \Delta_i^c$. Hence $x_1(t) \sim x_2(t)$ for all $t \notin \Delta$, because
$R$ is closed.

Observe next that if $t \in \Delta$, then there exists a sequence $t_k \in (t,\infty) \cap
\Delta^c$ such that $t_k \to t$. Then $x_i(t_k) \to x_i(t)$ by the right-continuity of $x_i$, and
again the fact that $R$ is closed implies $x_1(t) \sim x_2(t)$.
\end{proof}

\begin{theorem}
\label{the:fidiCT}
Two random processes $X$ and $Y$ with paths in $D_1(\R_+,S_1)$ and $D_2(\R_+,S_2)$, respectively,
satisfy $X \simst Y$ if and only if $(X_{t_1}, \dots, X_{t_n}) \simst (Y_{t_1}, \dots, Y_{t_n})$
for all finite parameter combinations $(t_1,\dots,t_n)$.
\end{theorem}
\begin{proof}
The necessity is obvious. To prove the converse, define for each positive integer $m$ the
discretization map $\pi_i^m: D_i \to S_i^{m^2+1}$ by
\[
 \pi_i^m(x) = (x(k/m))_{k=0}^{m^2},
\]
and the corresponding interpolation map $\eta_i^m: S_i^{m^2+1} \to D_i$ by
\[
 \eta_i^m(\alpha)(t) = \left\{
 \begin{aligned}
 \alpha_k,     &\quad t \in [k/m, (k+1)/m), \quad k = 0,1,\dots,m^2-1, \\
 \alpha_{m^2}, &\quad t \in (t,\infty).
 \end{aligned}
 \right.
\]
The functions $\pi_i^m$ and $\eta_i^m$ are measurable with respect to the Borel $\sigma$-algebras
on $D_i$ and $S_i^{m^2+1}$~\cite[Proposition 3.7.1]{ethier1986}. Let $(\hat X_1^m, \hat X_2^m)$ be
a coupling of $\pi_1^m(X_1)$ and $\pi_2^m(X_2)$ such that $(\hat X_1^m, \hat X_2^m) \in R^{m^2+1}$
almost surely. Then $(\eta_1^m(\hat X_1^m), \eta_2^m(\hat X_2^m))$ couples $\eta_1^m \circ
\pi_1^m(X_1)$ and $\eta_2^m \circ \pi_2^m(X_2)$, and moreover, $(\eta_1^m(\hat X_1^m),
\eta_2^m(\hat X_2^m)) \in R^D$ almost surely. Hence
\[
 \eta_1^m \circ \pi_1^m(X_1) \simst \eta_2^m \circ \pi_2^m(X_2).
\]
Because $\eta_i^m \circ \pi_i^m(x_i)$ converges to $x_i$ in $D_i$ as $m \to \infty$ for all $x_i
\in D_i$ (\cite[Problem 3.12]{ethier1986}, \cite[Lemma 3]{billingsley1999}), it follows that
$\eta_i^m \circ \pi_i^m(X_i) \wto X_i$. Lemmas~\ref{the:closed} and~\ref{the:closedSkorohod} show
that $\Rst^D$ is closed in the weak topology, so that $X_1 \simst X_2$.
\end{proof}

\subsection{Markov processes}
\label{sec:markov}

In the sequel, the notation $X(\mu,t)$ refers to the state of a Markov process $X$ at time $t$
with initial distribution $\mu$, and we shall use $X(x,t)$ as shorthand for $X(\delta_x,t)$.
Markov processes $X_1$ and $X_2$ are said to \emph{stochastically preserve} a relation $R$, if for
all $t$,
\[
 x_1 \sim x_2 \implies X_1(x_1,t) \simst X_2(x_2,t),
\]
or equivalently (see Theorem~\ref{the:preservation}),
\[
 \mu_1 \simst \mu_2 \implies X_1(\mu_1,t) \simst X_2(\mu_2,t).
\]
The following theorem presents a simple but powerful result, which together with the subrelation
algorithm (see Theorem~\ref{the:subrelationDT} below) provides a method for stochastically
relating (potentially unknown) stationary distributions of Markov processes based on their
generators.

\begin{theorem}
\label{the:stationaryDistributions} Let $X_1$ and $X_2$ be Markov processes with stationary
distributions $\mu_1$ and $\mu_2$ such that $X_i(x_i,t) \wto \mu_i$ as $t \to \infty$ for all
initial states $x_i$. Given any measurable relation $R$, a sufficient condition for $\mu_1 \simst
\mu_2$ is that $X_1$ and $X_2$ stochastically preserve some nontrivial closed subrelation of $R$.
\end{theorem}
\begin{proof}
Choose a pair of initial states $(x_1,x_2) \in R'$, where $R'$ is a closed subrelation of $R$
stochastically preserved by $X_1$ and $X_2$. Then $X_1(x_1,t)$ and $X_2(x_2,t)$ are stochastically
related with respect to $R'$ for all $t$, so by Lemma~\ref{the:closed} we see that $(\mu_1,\mu_2)
\in \Rst'$. Because $\Rst' \subset \Rst$, it follows that $(\mu_1,\mu_2) \in \Rst$.
\end{proof}

Let $X_1$ and $X_2$ be discrete-time Markov processes with transition probability kernels $P_1$
and $P_2$, respectively. The following result characterizes precisely when $X_1$ and $X_2$
stochastically preserve a relation $R$. A Markov process $\hat X$ taking values in $S_1 \times
S_2$ is called a \emph{Markovian coupling} of $X_1$ and $X_2$, if $\hat X(x,t)$ couples
$X_1(x_1,t)$ and $X_2(x_2,t)$ for all $t$ and all $x=(x_1,x_2)$. A measurable set $B$ is called
\emph{invariant} for a Markov process $X$, if $x \in B$ implies $X(x,t) \in B$ for all $t$ almost
surely.

\begin{theorem}
\label{the:DTRelation}
The following are equivalent:
\begin{enumerate}[(i)]
 \item \label{pro:markovDT1} $X_1$ and $X_2$ stochastically preserve the relation $R$.
 \item \label{pro:markovDT2} $P_1(x_1,B) \le P_2(x_2,B^\ra)$ for all $x_1 \sim x_2$ and compact $B \subset S_1$.
 \item \label{pro:markovDT3} $P_1$ and $P_2$ stochastically preserve the relation $R$.
 \item \label{pro:markovDT4} There is a Markovian coupling of $X_1$ and $X_2$ for which $R$ is invariant.
\end{enumerate}
\end{theorem}
\begin{proof}
The implications \reff{pro:markovDT4}$\implies$\reff{pro:markovDT1}$\implies$\reff{pro:markovDT3}
are direct consequences of the definitions, while \reff{pro:markovDT2} $\Longleftrightarrow$
\reff{pro:markovDT3} follows by Theorem~\ref{the:preservation}. For
\reff{pro:markovDT3}$\implies$\reff{pro:markovDT4}, observe that
Theorem~\ref{the:kernelCouplingGeneral} implies the existence of a coupling $P$ of the probability
kernels $P_1$ and $P_2$ such that $P(x,R) = 1$ for all $x \in R$. Let $\hat X$ be a discrete-time
Markov process with transition probability kernel $P$. Induction then shows that $\hat X$ is a
Markovian coupling of $X_1$ and $X_2$ for which $R$ is invariant.
\end{proof}

Theorems~\ref{the:subrelationKernel} and~\ref{the:DTRelation} yield the following characterization
for subrelations of a closed relation $R$ that are stochastically preserved by discrete-time
Markov processes $X_1$ and $X_2$ with continuous transition probability kernels $P_1$ and $P_2$.
Denote by $R^*$ the output~\reff{eq:RStar} of the subrelation algorithm in
Section~\ref{sec:subrelation}.

\begin{theorem}
\label{the:subrelationDT} $R^*$ is the maximal closed subrelation of $R$ that is stochastically
preserved by $X_1$ and $X_2$. Especially, $X_1$ and $X_2$ stochastically preserve a nontrivial
closed subrelation of $R$ if and only if $R^* \neq \emptyset$.
\end{theorem}

Markov jump processes shall be consider next. Recall that a map $Q: S \times \mB(S) \to \R_+$ is
called a \emph{rate kernel} on $S$, if $Q(x,dy) = q(x) P(x,dy)$ for some probability kernel $P$
and a positive measurable function $q$. A rate kernel $Q$ is called \emph{nonexplosive}, if the
standard construction using a discrete-time Markov process with transition probability kernel $P$
generates a Markov jump process with paths in $D(\R_+,S)$~\cite[Theorem 12.18]{kallenberg2002}.

Theorem~\ref{the:CTRelation} below characterizes precisely when a pair of Markov jump processes
$X_1$ and $X_2$ with nonexplosive rate kernels $Q_1$ and $Q_2$ stochastically preserves a closed
relation $R$. The construction of the Markovian coupling is based on the local uniformization of
the rate kernels $Q_i$ using the probability kernels $\hat P_i$ from $\Sprod$ to $S_i$, defined by
\begin{equation}
 \label{eq:hatPDef}
 \hat P_i(x, B_i)
 = \frac{q_i(x_i)}{q(x)} P_i(x_i, B_i) + \left(1 - \frac{q_i(x_i)}{q(x)} \right) \delta(x_i,B_i),
\end{equation}
where $q(x) = 1 + q_1(x_1) + q_2(x_2)$~\cite[Section 3]{lopez2000}.

\begin{theorem}
\label{the:CTRelation}
The following are equivalent:
\begin{enumerate}[(i)]
\item \label{pro:CTRelation1} $X_1$ and $X_2$ stochastically preserve the relation $R$.
\item \label{pro:CTRelation2}
For all $x_1 \sim x_2$ and compact $B \subset S_1$ such that $\delta(x_1,B) = \delta(x_2,B^\ra)$,
\[
 Q_1(x_1,B) - q_1(x_1) \delta(x_1,B) \le Q_2(x_2,B^\ra) - q_2(x_2) \delta(x_2,B^\ra).
\]
\item \label{pro:CTRelationNew} The probability kernels in~\reff{eq:hatPDef} satisfy $\hat P_1(x,\cdot) \simst \hat P_2(x,\cdot)$ for all $x \in R$.
\item \label{pro:CTRelation3} There is a Markovian coupling of $X_1$ and $X_2$ for which $R$ is invariant.
\end{enumerate}
\end{theorem}

Countable spaces admit the following slightly more convenient characterization, due to the fact
that all sets are measurable.

\begin{theorem}
\label{the:CTRelationCountable}
If the spaces $S_1$ and $S_2$ are countable, then the properties of Theorem~\ref{the:CTRelation}
are equivalent to requiring that for all $x_1 \sim x_2$:
\begin{equation}
 \label{eq:rateKernelRight}
 Q_1(x_1,B_1) \le Q_2(x_2,B_1^\ra)
\end{equation}
for all $B_1 \subset S_1$ such that $x_1 \notin B_1$ and $x_2 \notin B_1^\ra$, and
\begin{equation}
 \label{eq:rateKernelLeft}
 Q_1(x_1,B_2^\la) \ge Q_2(x_2,B_2)
\end{equation}
for all $B_2 \subset S_2$ such that $x_1 \notin B_2^\la$ and $x_2 \notin B_2$.
\end{theorem}
\begin{remark}
\label{rem:CTOrderCountable} For order relations on countable spaces it suffices to
verify~\reff{eq:rateKernelRight} for all upper sets $B_1$ and~\reff{eq:rateKernelLeft} for all
lower sets $B_2$ (see Remark~\ref{rem:stochasticOrder}), so Theorem~\ref{the:CTRelationCountable}
becomes equivalent to Massey's characterization~\cite[Theorem 3.4]{massey1987}.
\end{remark}

To prove Theorem~\ref{the:CTRelation} we need the following special form of
Theorem~\ref{the:kernelCouplingGeneral} to account for the fact that $\hat P_i$ are probability
kernels from $\Sprod$ to $S_i$, and not from $S_i$ to $S_i$.

\begin{lemma}
\label{the:kernelCouplingSpecial}
Let $P_i$ be probability kernels from $\Sprod$ to $S_i$ such that $P_1(x,\cdot) \simst
P_2(x,\cdot)$ for all $x \in R$. Then there exists a probability kernel $P$ on $\Sprod$ such that
$P(x,\cdot)$ couples $P_1(x,\cdot)$ and $P_2(x,\cdot)$ for all $x \in \Sprod$, and $P(x,R) = 1$
for all $x \in R$.
\end{lemma}
\begin{proof}
Denote $\hat S_i = \Sprod$ and $\hat S_i' = S_i$, $i=1,2$, and define $\hat R = \{(x,x) : x \in
R\}$ and $\hat R' = R$. Then by Theorem~\ref{the:kernelCouplingGeneral} there exists a probability
kernel $\hat P$ from $\hat S_1 \times \hat S_2$ to $\hat S_1' \times \hat S_2'$ such that $\hat
P((x,y), \cdot)$ couples $P_1(x,\cdot)$ and $P_2(y,\cdot)$ for all $x \in \Sprod$ and $y \in
\Sprod$, and $\hat P( (x,x), R) = 1$ for all $x \in R$. Define $P(x,B) = \hat P((x,x),B)$ for all
$x \in \Sprod$ and measurable $B \subset \Sprod$.
\end{proof}

\begin{proof}[Proof of Theorem~\ref{the:CTRelation}]
\reff{pro:CTRelation1} $\implies$ \reff{pro:CTRelation2}. Choose $x_1 \sim x_2$ and a compact $B
\subset S_1$ such that $\delta(x_1,B) = \delta(x_2,B^\ra)$. Then Theorem~\ref{the:relationTest}
shows that
\[
 \pr(X_1(x_1,t) \in B) - \delta(x_1,B) \le \pr(X_2(x_2,t) \in B^\ra) - \delta(x_2,B^\ra)
\]
for all $t$. Dividing both sides above by $t$, and taking $t \downarrow 0$, we thus see using
Kolmogorov's backward equation~\cite[Theorem 12.25]{kallenberg2002} the validity
of~\reff{pro:CTRelation2}.

\reff{pro:CTRelation2} $\implies$ \reff{pro:CTRelationNew}. Observe that for any compact $B
\subset S_1$,
\begin{align*}
 \hat P_2(x,B^\ra) - \hat P_1(x,B) &=  q(x)^{-1}  \{ ( Q_2(x_2,B^\ra) - q(x_2)\delta(x_2,B^\ra) ) \\
 &\qquad -  (Q_1(x_1,B) - q(x_1)\delta(x_1,B)) \} \\
 &\qquad + \delta(x_2,B^\ra) - \delta(x_1,B).
\end{align*}
Because $\delta(x_2,B^\ra) \ge \delta(x_1,B)$ for all $x_1 \sim x_2$, we see that $\hat P_1(x,B)
\le \hat P_2(x,B^\ra)$ for all $x \in R$. Hence using Theorem~\ref{the:relationTest} we may
conclude that $\hat P_1(x,\cdot) \simst \hat P_2(x,\cdot)$ for all $x \in R$.

\reff{pro:CTRelationNew} $\implies$ \reff{pro:CTRelation3}. Let $\hat P_i$ be the probability
kernels defined in~\reff{eq:hatPDef}. Lemma~\ref{the:kernelCouplingSpecial} then shows the
existence of a probability kernel $P$ on $\Sprod$ such that $P(x,\cdot)$ couples $\hat
P_1(x,\cdot)$ and $\hat P_2(x,\cdot)$ for all $x \in \Sprod$, and $P(x,R) = 1$ for all $x \in R$.
Define a rate kernel $Q$ on $\Sprod$ by $ Q(x,B) = q(x) P(x,B)$. Then
\[
 \int_{\Sprod} \left( f_i(y_i) - f_i(x_i) \right) \, Q(x,dy)
 = \int_{S_i} \left( f_i(y_i) - f_i(x_i) \right) \, Q_i(x_i,dy_i)
\]
for all $x \in \Sprod$ and all bounded measurable $f$ on $S_i$, $i=1,2$, which implies that $Q$ is
nonexplosive (Chen~\cite[Theorem 37]{chen1986}). Let $\hat X$ be a Markov jump process generated
by $Q$ using the standard construction~\cite[Theorem 12.18]{kallenberg2002}. Then $\hat X(x,t)$
couples $X_1(x_1,t)$ and $X_2(x_2,t)$ for all $x \in \Sprod$ and for all $t$ (Chen~\cite[Theorem
13]{chen1986}). Moreover, $R$ is invariant for $\hat X$, because $P(x,R) = 1$ for all $x \in R$.

\reff{pro:CTRelation3}$\implies$\reff{pro:CTRelation1}. Clear by definition.
\end{proof}

\begin{proof}[Proof of Theorem~\ref{the:CTRelationCountable}]
Assume first that~\reff{eq:rateKernelRight} and~\reff{eq:rateKernelLeft} hold, and choose $x_1
\sim x_2$ and $B_1 \subset S_1$ so that $\delta(x_1,B_1) = \delta(x_2,B_1^\ra)$. If $x_1 \notin
B_1$, and $x_2 \notin B_1^\ra$, then the validity of
Theorem~\ref{the:CTRelation}:\reff{pro:CTRelation2} is obvious from~\reff{eq:rateKernelRight}. In
the other case where $x_1 \in B_1$, and $x_2 \in B_1^\ra$, denote $B_2 = (B_1^\ra)^c$. Then it
follows (check) that $B_2^\la \subset B_1^c$, and thus $x_1 \notin B_2^\la$. This further implies
that $x_2 \notin B_2$, because $x_1 \sim x_2$. Hence~\reff{eq:rateKernelLeft} shows that
\[
 Q_2(x_2,B_2)
 \le Q_1(x_1,B_2^\la)
 \le Q_1(x_1,B_1^c),
\]
from which we again see that Theorem~\ref{the:CTRelation}:\reff{pro:CTRelation2} is valid.

Assume next that Theorem~\ref{the:CTRelation}:\reff{pro:CTRelation3} holds, and let $x_1 \sim x_2$
be such that $x_1 \notin B_1$ and $x_2 \notin B_1^\ra$. Then~\reff{eq:rateKernelRight} follows by
Theorem~\ref{the:relationTest}, because $Q_1(x_1,B_1) = \hat P_1(x,B_1)$ and $Q_2(x_2,B_1^\ra) =
\hat P_2(x_2,B_1^\ra)$. Inequality~\reff{eq:rateKernelLeft} can be verified by a symmetrical
argument.
\end{proof}

This section is concluded by an analogue of Theorem~\ref{the:subrelationDT}. A rate kernel
$Q(x,dy) = q(x) P(x,dy)$ such that $q$ is a continuous function and $P$ is a continuous
probability kernel shall be called \emph{continuous}. Given Markov jump processes $X_1$ and $X_2$
with continuous nonexplosive rate kernels $Q_1$ and $Q_2$, define the relations $R^{(n)}$ by
$R^{(0)} = R$, and
\begin{equation}
\label{eq:subrelationCTKernel}
R^{(n+1)} = \left\{x \in R^{(n)}: \ (\hat P_1(x,\cdot), \hat P_2(x,\cdot)) \in \Rst^{(n)}
\right\},
\end{equation}
where $\hat P_1$ and $\hat P_2$ are given by~\reff{eq:hatPDef}. Moreover, denote $R^*=
\cap_{n=0}^\infty R^{(n)}$, as in Section~\ref{sec:subrelation}.

\begin{theorem}
\label{the:subrelationCT} $R^*$ is the maximal closed subrelation of $R$ that is stochastically
preserved by $X_1$ and $X_2$. Especially, $X_1$ and $X_2$ stochastically preserve a nontrivial
closed subrelation of $R$ if and only if $R^* \neq \emptyset$.
\end{theorem}
\begin{proof}
The continuity of $Q_1$ and $Q_2$ guarantees the continuity of $\hat P_1$ and $\hat P_2$. The
proof is hence completed by repeating the steps in the proof of
Theorem~\ref{the:subrelationKernel}, with notational modifications to take into account that here
$\hat P_i$ are probability kernels from $\Sprod$ to $S_i$, and not from $S_i$ to $S_i$.
\end{proof}

\section{Applications}
\label{sec:applications}

\subsection{Hidden Markov processes}

The goal of this section is to stochastically compare two hidden Markov processes $Y_1$ and $Y_2$
of the form $Y_i = f_i \circ X_i$, where $X_i$ is a Markov process taking values in $S_i'$, and
$f_i$ is a continuous function from $S_i'$ to $S_i$. Although the results in this section have
natural counterparts for Markov jump processes, we shall only treat the case where $X_i$ are
discrete-time Markov processes with transition probability kernels $P_i$.

\begin{theorem}
\label{the:inducedKernelTest}
Let $R$ be a closed relation between $S_1$ and $S_2$, and assume that
\begin{equation}
 \label{eq:inducedKernelTest}
 P_1(x_1, f_1^{-1}(B)) \le P_2(x_2, f_2^{-1}(B^\ra))
\end{equation}
for all $x_1$ and $x_2$ such that $f_1(x_1) \sim f_2(x_2)$ and all compact $B \subset S_1$. Then
$Y_1(t) \simst Y_2(t)$ whenever $Y_1(0) \simst Y_2(0)$.
\end{theorem}

The proof of Theorem~\ref{the:inducedKernelTest} relies on the notion of induced relation, defined
as follows. Given a closed relation $R$ between $S_1$ and $S_2$, and two continuous functions
$\phi_i: S_i' \to S_i$, $i=1,2$, we denote
\[
 R' = \left\{ (x_1',x_2') \in \Sprodp: (\phi_1(x_1'), \phi_2(x_2')) \in R \right\}.
\]
The closed relation $R'$ is said to be \emph{induced from $R$} by the pair $(\phi_1,\phi_2)$, and
we denote by $\Rst'$ the stochastic relation generated by $R'$.

\begin{lemma}
\label{the:inducedRelationNew} Probability measures $\mu_1$ on $S_1'$ and $\mu_2$ on $S_2'$
satisfy $(\mu_1,\mu_2) \in \Rst'$ if and only if $(\mu_1 \circ \phi_1^{-1}, \mu_2 \circ
\phi_2^{-1}) \in \Rst$.
\end{lemma}
\begin{proof}
Observe that $R' = \phi^{-1}(R)$, where the function $\phi$ is defined by $\phi(x_1,x_2) =
(\phi_1(x_1),\phi_2(x_2))$. If $(\mu_1,\mu_2) \in \Rst'$, then let $\lambda$ be a coupling of
$\mu_1$ and $\mu_2$ such that $\lambda(R') = 1$. Then $\lambda \circ \phi^{-1}$ couples the
probability measures $\mu_1 \circ \phi_1^{-1}$ and $\mu_2 \circ \phi_2^{-1}$, and moreover, that
$\lambda \circ \phi^{-1}(R) = 1$.

Assume next that $(\nu_1,\nu_2) \in \Rst$, where $\nu_i = \mu_i \circ \phi_i^{-1}$. Choose a
compact $B' \subset S_1'$, and observe that the right conjugate of $B'$ with respect to $R'$
equals
\[
 [B']_{R'}^\ra = \phi_2^{-1}( [\phi_1(B')]_{R}^\ra),
\]
where $[\phi_1(B')]_{R}^\ra$ denotes the right conjugate of $\phi_1(B')$ with respect to $R$.
Moreover, because $\phi_1(B')$ is compact, Theorem~\ref{the:relationTest} shows that
\[
 \mu_1(B') = \nu_1(\phi_1(B')) \le \nu_2([\phi_1(B')]_{R}^\ra) = \mu_2([B']_{R'}^\ra),
\]
and hence $(\mu_1, \mu_2) \in \Rst'$.
\end{proof}

\begin{proof}[Proof of Theorem~\ref{the:inducedKernelTest}]
Let $R'$ be the relation between $S_1'$ and $S_2'$ induced from $R$ by the pair $(f_1,f_2)$.
Lemma~\ref{the:inducedRelationNew} then shows that for all $t$,
\[
 Y_1(t) \stochrel{R} Y_2(t) \quad\text{if and only if} \quad
 X_1(t) \stochrel{R'}  X_2(t),
\]
and moreover, \reff{eq:inducedKernelTest} is equivalent to requiring that the pair $(P_1,P_2)$
stochastically preserves $R'$.
\end{proof}

\begin{example}[Non-Markov processes]
Let $Y_1$ and $Y_2$ be non-Markov processes with values in $S_1$ and $S_2$, respectively.
Following Whitt~\cite{whitt1986}, let us assume that $Y_i$ can be made Markov by keeping track of
additional information, say a random processes $Z_i$ with values in $S_i'$.  Under this
assumption, $Y_i = \pi_i \circ X_i$, where $X_i = (Y_i,Z_i)$ is a Markov process taking values
$S_i \times S_i'$, and $\pi_i$ denotes the projection from $S_i \times S_i'$ onto $S_i$. Denoting
the transition probability kernel of $X_i$ by $P_i$, condition~\reff{eq:inducedKernelTest} in
Theorem \ref{the:inducedKernelTest} becomes equivalent to
\begin{equation}
\label{eq:whitt} \sup_{z_1 \in S_1'} P_1((y_1,z_1), \, B \times S_1') \le \inf_{z_2 \in S_2'}
P_2((y_2,z_2), \, B^\ra \! \times S_2')
\end{equation}
for all $y_1 \sim y_2$ and all compact $B \subset S_1$. Inequality~\reff{eq:whitt} together with
$Y_1(0) \simst Y_2(0)$ is thus sufficient for $Y_1(t) \simst Y_2(t)$ for all $t$. This formulation
is conceptually similar, though not equivalent, to~\cite[Theorem 1]{whitt1986}.
\end{example}

\begin{example}[Lumpability]
A Markov process $X$ with values in $S_1$ is called \emph{lumpable} with respect to $f:S_1 \to
S_2$, if $f \circ X$ is Markov for any initial distribution of $X$ \cite{kemeny1960}. It is well
known \cite{dynkin1965} that $X$ is lumpable if and only if its transition probability kernel
satisfies
\[
 P(x,f^{-1}(B)) = P(y,f^{-1}(B))
\]
for all $x$ and $y$ such that $f(x) = f(y)$ and all measurable $B \subset S_2$. This is equivalent
to saying that the pair $(P,P)$ stochastically preserves the relation $\{(x,y) \in S_1 \times S_1:
f(x) = f(y)\}$ induced by $(f,f)$ from the equality on $S_2$. Moreover, if $X$ is lumpable with
respect to $f$, then the pair $(P,P')$ stochastically preserves the relation $\{(x,y) \in S_1
\times S_2: f(x) = y\}$ induced from the equality on $S_2$ by the pair $(f,\idmap)$, where $P'$
denotes the transition probability kernel of $f \circ X$, and $\idmap$ is the identity map on
$S_2$. The notion of lumpability may be generalized by calling $X$ lumpable with respect to a
relation $R$, if there exists a Markov process $Y$ such that $X$ and $Y$ stochastically preserve
$R$~\cite{lopez2002}.
\end{example}

\begin{example}[Stochastic induced order]
For order relations, condition~\reff{eq:inducedKernelTest} in Theorem~\ref{the:inducedKernelTest}
can be rephrased as
\[
 P_1(x_1, f_1^{-1}(B)) \le P_2(x_2, f_2^{-1}(B))
\]
for all upper sets $B \subset S$ and all $x_1$ and $x_2$ such that $f_1(x_1) \le f_2(x_2)$ (see
Remark~\ref{rem:stochasticOrder}). This is a discrete-time analogue to~\cite[Theorem
6]{doisy2000}.
\end{example}

\subsection{Population processes}
\label{sec:populationProcesses}

Let us denote by $e_1,\dots,e_m$ the unit vectors of $\Z^m$, and define $e_0 = 0$ and $e_{i,j} =
-e_i + e_j$ for notational convenience. A \emph{Markov population process}~\cite{kingman1969} is a
nonexplosive Markov jump process taking values\footnote{Often it is natural to assume that the
elements of $S$ have positive coordinates, but in this section there is no need to make this
restriction.} in $S \subset \Z^m$, generated by the transitions
\[
 x \mapsto x + e_{i,j} \quad \text{at rate} \ \alpha_{i,j}(x), \quad i,j \in \{0,1,\dots,m\}, \quad i \neq
 j.
\]
where $\alpha_{i,j}$ are positive functions on $S$ such that $\alpha_{i,j}(x) = 0$ for $x +
e_{i,j} \notin S$. The functions $\alpha_{0,i}$ and $\alpha_{i,0}$ may be regarded as the arrival
and departure rates of individuals for colony $i$, and $\alpha_{i,j}$ represents the transfer rate
of individuals from colony $i$ to colony $j$. Population processes with values in $S \subset
\Z_+^m$ may be viewed as Markovian queueing networks~\cite{serfozo1999} or interacting particle
systems~\cite{lindvall1992}.

The following result characterizes precisely when two population processes stochastically preserve
a relation $R$ between $S \subset \Z^m$ and $S' \subset \Z^{m'}$. To state the result, define for
$x \in S$ and $y \in S'$, and for sets of index pairs $U \subset \{0,\dots,m\}^2$ and $V \subset
\{0,\dots,m'\}^2$,
\[
 U_\rightarrow(x,y) = \left\{ (k,l): x + e_{i,j} \sim y + e_{k,l} \ \text{for some} \ (i,j) \in U \right\}
\]
and
\[
 V_\leftarrow(x,y) = \left\{ (i,j): x + e_{i,j} \sim y + e_{k,l} \ \text{for some} \ (k,l) \in V \right\}.
\]

\begin{theorem}
\label{the:populationRelation}
Let $X$ and $X'$ be population processes taking values in $S \subset \Z^m$ and $S' \subset
\Z^{m'}$ generated by transition rate functions $\alpha_{i,j}$ and $\alpha_{k,l}'$, respectively.
Then $X$ and $X'$ stochastically preserve a relation $R$ if and only if for all $x \sim y$:
\begin{equation}
 \label{eq:populationRelationRight}
 \sum_{(i,j) \in U} \alpha_{i,j}(x) \le \sum_{(k,l) \in U_\rightarrow(x,y)} \alpha_{k,l}'(y)
\end{equation}
for all $U \subset \{(i,j): \ x + e_{i,j} \not\sim y\}$, and
\begin{equation}
 \label{eq:populationRelationLeft}
 \sum_{(i,j) \in V_\leftarrow(x,y)} \alpha_{i,j}(x) \ge \sum_{(k,l) \in V} \alpha_{k,l}'(y)
\end{equation}
for all $V \subset \{(k,l): \ x \not\sim y + e_{k,l}\}$.
\end{theorem}
\begin{proof}
Observe that for any $x \sim y$ and $B \subset S$ such that $x \notin B$ and $y \notin B^\ra$, the
rate kernels of $X$ and $X'$ satisfy
\[
 Q(x,B)      = \sum_{(i,j) \in U} \alpha_{i,j}(x)
 \quad \text{and} \quad
 Q'(y,B^\ra) = \sum_{(k,l) \in U_\ra(x,y)} \alpha'_{k,l}(y),
\]
where $U = \{(i,j): x + e_{i,j} \in B\}$. Moreover, $x + e_{i,j} \not\sim y$ for all $(i,j) \in
U$, because $y \notin B^\ra$. Hence inequality~\reff{eq:populationRelationRight} is equivalent
to~\reff{eq:rateKernelRight} in Theorem~\ref{the:CTRelationCountable}. By symmetry, we see
that~\reff{eq:populationRelationLeft} is equivalent~\reff{eq:rateKernelLeft}.
\end{proof}

\begin{example}[Partial coordinatewise order]
\label{exa:randomWalks}
Define a relation between $S$ and $S'$ by denoting $x \leM y$, if $x_i \le y_i$ for all $i \in M$,
where $M$ is a fixed subset of $\{1,\dots,m\} \cap \{1,\dots,m'\}$. We shall next show that $X$
and $X'$ stochastically preserve the relation $\leM$ if and only if for all $x \leM y$ and all $k$
in the set $M_0(x,y) = \{i \in M: x_i = y_i\}$:
\begin{equation}
 \label{eq:subpopulationUpper}
 \sum_{i \in I} \alpha_{i,k}(x)
 \le
 \sum_{i \in I \cup ([0,m'] \setminus M_0(x,y))} \alpha'_{i,k}(y)
\end{equation}
for all $I \subset [0,m] \setminus \{k\}$, and
\begin{equation}
 \label{eq:subpopulationLower}
 \sum_{j \in J \cup ([0,m] \setminus M_0(x,y))} \alpha_{k,j}(x)
 \ge \sum_{j \in J} \alpha'_{k,j}(y)
\end{equation}
for all $J \subset [0,m']\setminus \{k\}$.

To justify the above claim, observe that for any $x \leM y$ and $U$ is as in
Theorem~\ref{the:populationRelation}, we can write $U = \cup_{k \in K} (I_k \times \{k\})$ for
some $I_k \subset [0,M] \setminus \{k\}$ and $K \subset M_0(x,y)$. Hence $U_\ra(x,y) = \cup_{k \in
K} (J_k \times \{k\})$, where $U_\ra(x,y) = \cup_{k \in K} (I_k \cup ([0,m'] \setminus M_0(x,y)))
\times \{k\})$. By summing both sides of~\reff{eq:subpopulationUpper} over $k \in K$, with $I_k$
in place of $I$, we see that~\reff{eq:populationRelationRight} holds. On the other hand, it is
easy to see that~\reff{eq:populationRelationRight} implies~\reff{eq:subpopulationUpper}, and a
symmetric argument shows the equivalence of~\reff{eq:populationRelationLeft}
and~\reff{eq:subpopulationLower}.

The above characterization extends and sharpens earlier comparison results for Markovian queueing
networks~\cite{economou2003, lindvall1992,lopez2000}. For population processes where transfers
between colonies do not occur, \reff{eq:subpopulationUpper} and~\reff{eq:subpopulationLower}
simplify to
\[
 \alpha_{0,k}(x) \le \alpha'_{0,k}(y)
 \quad\text{and}\quad
 \alpha_{k,0}(x) \ge \alpha'_{k,0}(y)
\]
for all $x \leM y$ and all $k \in M$ such that $x_k = y_k$~\cite[Lemma 1]{borstNow}.
\end{example}

\subsection{Parallel queueing system}
\label{sec:parallelQueueingSystem} Consider a system of two queues in parallel, where customers
arrive to queue $k$ at rate $\lambda_k \in (0,1)$ and have unit service rate. Assuming all
interarrival and service times are exponential, the queue length process $X=(X_1,X_2)$ is a Markov
population process on $\Z_+^2$ with transition rates $\alpha_{0,k}(x) = \lambda_k$ and
$\alpha_{k,0}(x) = 1(x_k > 0)$, $k=1,2$. We shall also consider a modification of the system,
where load is balanced by routing incoming traffic to the shortest queue, modeled as a Markov
population process $\XLB = (\XLB_1,\XLB_2)$ with transition rates
\begin{align*}
 \alphaLB_{0,1}(x) &= (\lambda_1 + \lambda_2) 1(x_1 < x_2) + \lambda_1 1(x_1 = x_2), \\
 \alphaLB_{0,2}(x) &= (\lambda_1 + \lambda_2) 1(x_1 > x_2) + \lambda_2 1(x_1 = x_2),
\end{align*}
and $\alphaLB_{k,0}(x) = \alpha_{k,0}(x)$ for $k=1,2$. Common sense suggests that load balancing
decreases the total number of customers in the system, so that
\begin{equation}
 \label{eq:loadBalancing}
 \XLB_1(t) + \XLB_2(t) \lest X_1(t) + X_2(t).
\end{equation}
However, the justification of~\reff{eq:loadBalancing} appears difficult, because using
Theorem~\ref{the:populationRelation} can check that the processes $\XLB$ and $X$ do not
stochastically preserve the coordinatewise order on $\Z_+^2$, nor the order
\[
\Rsum = \{(x,y): |x| \le |y|\},
\]
where $|x| = x_1 + x_2$. On the other hand, it is known~\cite{winston1977}
that~\reff{eq:loadBalancing} holds for all $t$, whenever $\XLB(0) = X(0)$, which suggests that
$\XLB$ and $X$ might stochastically preserve some strict subrelation of $\Rsum$. The following
theorem summarizes the output of the subrelation algorithm applied to the rate kernels of $\XLB$
and $X$.

\begin{theorem}
\label{the:queueing} Starting from $R^{(0)} = \Rsum$, the subrelation
iteration~\reff{eq:subrelationCTKernel} produces the sequence of relations
\begin{equation}
 \label{eq:queueing}
 R^{(n)} = \left\{(x,y): |x| \le |y| \ \text{and} \ x_1 \vee x_2 \le y_1 \vee y_2 + (y_1 \wedge y_2 - n)^+
 \right\},
\end{equation}
which converges to
\[
 R^* = \left\{(x,y): |x| \le |y| \ \text{and} \ x_1 \vee x_2 \le y_1 \vee y_2 \right\}.
\]
\end{theorem}

The limiting relation $R^*$ may be identified as the weak majorization order $\wmaj$ on $\Z_+^2$
(Example~\ref{exa:majorization}). As a consequence,
\[
 \XLB(0) \wmaj X(0)
 \implies
 \XLB(t) \wmajst X(t) \quad \text{for all $t$}.
\]
Especially, $\XLB(0) \wmaj X(0)$ implies~\reff{eq:loadBalancing}, and moreover,
\[
 \XLB_1(t) \vee \XLB_2(t) \, \lest \, X_1(t) \vee X_2(t),
\]
which indicates that the queue lengths corresponding to $\XLB$ are more balanced than those
corresponding to $X$. The proof of Theorem~\ref{the:queueing} is based on the following lemma, the
proof of which is omitted.
\begin{lemma}
\label{the:alpha} The function $\alpha_n(x) = x_1 \vee x_2 \vee (|x|-n)$ on $\Z^2$ satisfies:
\begin{enumerate}[(i)]
\item \label{pro:alphaEq} $\alpha_n(x) = |x| - x_1 \wedge x_2 \wedge n$.
\item \label{pro:alphaLe} $\alpha_n(x-e_k) \le \alpha_{n+1}(x)$ for all $k$.
\item \label{pro:alphaAsymHi} $\alpha_n(x+e_1) > \alpha_n(x+e_2)$ if and only if $x_1 > x_2 \vee (|x|-n)$.
\item \label{pro:alphaAsymLo} $\alpha_n(x-e_1) < \alpha_n(x-e_2)$ if and only if $x_1 > x_2 \vee (|x|-n-1)$.
\item \label{pro:alphaVee} $\alpha_n(x+e_k) > \alpha_n(x)$ for $k$ such that $x_k = x_1 \vee x_2$.
\item \label{pro:alphaWedge} $\alpha_n(x-e_k) = \alpha_{n+1}(x)$ for $k$ such that $x_k = x_1 \wedge x_2$.
\end{enumerate}
\end{lemma}

\begin{proof}[Proof of Theorem~\ref{the:queueing}]
Define the function $\alpha_n(x)$ as in Lemma~\ref{the:alpha}. Then the relations $R^{(n)}$
defined in~\reff{eq:queueing} can be written as
\[
R^{(n)} = \left\{(x,y): |x| \le |y|, \ \alpha_n(x) \le \alpha_n(y) \right\},
\]
which shows that $R^{(n)}$ is an order for all $n$. To show that $R^{(n)}$ is the sequence of
relations produced by iteration~\reff{eq:subrelationCTKernel}, it is sufficient
(Remark~\ref{rem:CTOrderCountable}) to show that the properties
\begin{itemize}
\item $\QLB(x,U) \le Q(y,U)$ for all $R^{(n)}$-upper sets $U$ such that $x,y \notin U$,
\item $\QLB(x,U) \ge Q(y,U)$ for all $R^{(n)}$-lower sets $U$ such that $x,y \notin U$,
\end{itemize}
are valid for all $(x,y) \in R^{(n+1)}$, and that at least one of the above properties fails for
$(x,y) \in R^{(n)} \setminus R^{(n+1)}$.

Assume $(x,y) \in R^{(n+1)}$ (actually, $R^{(n)}$ is here enough), and let $U$ be an
$R^{(n)}$-upper set such that $x,y \notin U$. Then $x-e_k, y-e_k \notin U$ for all $k$.
\begin{enumerate}[(i)]

\item If $y+e_1, y+e_2 \in U$, then $\QLB(x,U) \le \lambda_1 + \lambda_2 = Q(y,U)$.

\item Assume $y+e_1, y+e_2 \notin U$, and choose $l$ so that $y_l = y_1 \vee y_2$. Then
Lemma~\ref{the:alpha}:\reff{pro:alphaVee} together with $\alpha_n(x) \le \alpha_n(y)$ imply that
$\alpha_n(x+e_k) \le \alpha_n(x) + 1 \le \alpha_n(y + e_l)$ for all $k$. Hence $x + e_k \notin U$
for all $k$, because $U$ is $R^{(n)}$-upper, so that $Q(x,U) = 0 = \QLB(y,U)$.

\item Assume $y+e_1 \in U, y+e_2 \notin U$. Then
$\alpha_n(y+e_1)> \alpha_n(y+e_2)$, so Lemma~\ref{the:alpha}:\reff{pro:alphaAsymHi} shows that
$y_1 > y_2$ and $y_1 > |y|-n$. Now if $x+e_k \in U$, then $\alpha_n(x+e_k) > \alpha_n(y+e_2) =
\alpha_n(y)$. Further $\alpha_n(x) \le \alpha_n(y)$ shows that $\alpha_n(x) = \alpha_n(y)$, so
that $\alpha_n(x) = y_1
> |x| - n$. Because $\alpha_n(x+e_k) > \alpha_n(x) > |x|-n$, it follows that $x_k = x_1 \vee x_2 = y_1$.
As a consequence, $|x| \le |y|$ implies that $x_1 \wedge x_2 \le y_2 < y_1$, so that $x_k > x_1
\wedge x_2$. Hence $\alphaLB_k(x) = 0$, which implies that $\QLB(x,U) = 0 \le Q(y,U)$. The case
where $y+e_1 \notin U, y+e_1 \in U$ is similar.
\end{enumerate}

Assume next that $(x,y) \in R^{(n+1)}$, and let $U$ be an $R^{(n)}$-lower set such that $x,y
\notin U$. Then $x+e_k, y+e_k \notin U$ for all $k$.
\begin{enumerate}[(i)]
\item Assume $y-e_1 \in U$, $y-e_2 \notin U$. Then $y_2 = 0$ or $\alpha_n(y-e_1) < \alpha_n(y-e_2)$,
so with the help of Lemma~\ref{the:alpha}:\reff{pro:alphaAsymLo}, we see that $y_1
> y_2 \vee (|y|-n-1)$. Choose $k$ so that $x_k = x_1 \vee x_2$, and observe that

$(x,y-e_1) \notin R^{(n)}$ implies that either $|x| = |y|$ or $\alpha_n(x) \ge \alpha_n(y)$, so
that $x_k > 0$. If $x_1 \neq x_2$, then $\alpha_n(x-e_k) = \alpha_n(x) - 1 \le \alpha_n(y) - 1 =
\alpha_n(y-e_1)$. If $x_1 = x_2$, then $|x| \le |y|$ together with $y_1 > y_2$ implies that $x_1 <
y_1$. Hence $\alpha_n(x-e_k) = x_1 \vee (|x|-n-1) \le (y_1-1) \vee (|y|-n-1) = \alpha_n(y-e_1)$.
Thus, we may conclude that $x-e_k \in U$, which shows that $\QLB(x,U) \ge 1 = Q(y,U)$. By
symmetry, the same conclusion holds under the assumption $y-e_1 \notin U$, $y-e_2 \in U$.

\item Assume $y-e_1, y-e_2 \in U$, and choose $l$ so that $y_l = y_1 \wedge y_2$.
Then using Lemma~\ref{the:alpha}:\reff{pro:alphaLe} and
Lemma~\ref{the:alpha}:\reff{pro:alphaWedge} we find that for all $k$,
\[
 \alpha_n(x-e_k) \le \alpha_{n+1}(x) \le \alpha_{n+1}(y) = \alpha_n(y-e_l).
\]
Further, if $|x| = |y|$, then Lemma~\ref{the:alpha}:\reff{pro:alphaEq} together with $x_1 \vee x_2
\le \alpha_{n+1}(x) \le \alpha_{n+1}(y)$ shows that $x_1 \wedge x_2 = |x| - x_1 \vee x_2 \ge y_1
\wedge y_2 \wedge (n+1)$, so that $x_1 \wedge x_2 \ge 1$. On the other hand, if $|x| < |y|$, then
$x \notin U$ implies that $\alpha_n(x) > \alpha_n(y-e_l) \ge \alpha_{n+1}(x)$, so it follows that
$x_1 \vee x_2 < |x| - n$, which again shows that $x_1 \wedge x_2 \ge 1$. Hence $x-e_k \in U$ for
all $k$, and we may conclude that $\QLB(x,U) = 2 = Q(y,U)$.

\item Assume $y-e_l \notin U$ for all $l$. Then $\QLB(x,U) \ge 0 = Q(y,U)$.
\end{enumerate}

Finally, assume that $(x,y) \in R^{(n)} \setminus R^{(n+1)}$. Then
\[
 \alpha_{n+1}(y) < \alpha_{n+1}(x)
 \le \alpha_n(x) \le \alpha_n(y),
\]
which implies that $y_1 \vee y_2 < |y|-n$ and $\alpha_{n}(x) = \alpha_n(y) = |y| - n$. Consider
the $R^{(n)}$-lower set $U = \{z: |z| < |y|, \ \alpha_n(z) < \alpha_n(y)\}$. Then $x,y \notin U$,
$x+e_k, y+e_k \notin U$ for all $k$, and $y-e_l \in U$ for all $l$. On the other hand,
$\alpha_n(x) = |y| - n$ together with $|x| \le |y|$ shows that $\alpha_n(x) = x_1 \vee x_2$, so
that $\alpha_n(x-e_k) = \alpha_n(x)$ for some $k$. Especially, $x-e_k \notin U$, so that
$\QLB(x,U) \le 1 < 2 = Q(y,U)$.
\end{proof}

\section{Conclusion}
\label{sec:conclusion}

This paper presented a systematic study of stochastic relations, which naturally extend the notion
of stochastic orders to relations between random variables and processes that may take values in
different state spaces. The key points of the paper may be summarized by
Theorem~\ref{the:subrelationKernel}, which characterizes the existence of subrelations
stochastically preserved by a pair of probability kernels, and
Theorem~\ref{the:stationaryDistributions}, which underlines the relevance of subrelation
techniques in stochastically comparing stationary distributions of Markov processes. Finite-state
Markov processes and diffusions are two important classes of processes that were not discussed in
the paper. In finite state spaces the subrelation algorithm converges in finite time, which calls
for numerical analysis of the runtime. The analysis of stochastic relations for diffusion
processes requires the identification of suitable test functions that behave well with respect to
taking relational conjugates. These issues may be considered interesting topics for future
research.

\section*{Acknowledgements}
The main part of the work presented has been carried out at Centrum voor Wiskunde en Informatica
(CWI) and Eindhoven University of Technology, the Netherlands. The research has been supported by
the Dutch BSIK/BRICKS PDC2.1 project, Helsingin Sanomat Foundation, and the Academy of Finland.

\bibliographystyle{apalike}
\bibliography{lslReferences-2008-06-19}

\newcommand{\SortNoop}[1]{}
\begin{thebibliography}{}

\bibitem[Billingsley, 1999]{billingsley1999}
Billingsley, P. (1999).
\newblock {\em Convergence of Probability Measures}.
\newblock Wiley, second edition.

\bibitem[Borst et~al., 2007]{borstNow}
Borst, S., Jonckheere, M., and Leskel\"a, L. (2007).
\newblock Stability of parallel queueing systems with coupled service rates.
\newblock Discrete Event Dyn. Syst., To appear.

\bibitem[Brandt and Last, 1994]{brandt1994}
Brandt, A. and Last, G. (1994).
\newblock On the pathwise comparison of jump processes driven by stochastic
  intensities.
\newblock {\em Math. Nachr.}, 167:21--42.

\bibitem[Chen, 1986]{chen1986}
Chen, M.-F. (1986).
\newblock Coupling for jump processes.
\newblock {\em Acta Math. Sin.}, 2(2):123--136.

\bibitem[Doisy, 2000]{doisy2000}
Doisy, M. (2000).
\newblock A coupling technique for stochastic comparison of functions of
  {Markov} processes.
\newblock {\em J. Appl. Math. Decis. Sci.}, 4(1):39--64.

\bibitem[Dynkin, 1965]{dynkin1965}
Dynkin, E.~B. (1965).
\newblock {\em Markov processes {Vol.} {I}}.
\newblock Springer.

\bibitem[Economou, 2003]{economou2003}
Economou, A. (2003).
\newblock Necessary and sufficient conditions for the stochastic comparison of
  {Jackson} networks.
\newblock {\em Probab. Eng. Inform. Sci.}, 17:143--151.

\bibitem[Eifler, 1976]{eifler1976}
Eifler, L.~Q. (1976).
\newblock Open mapping theorems for probability measures on metric spaces.
\newblock {\em Pac. J. Math.}, 66:89--97.

\bibitem[Ethier and Kurtz, 1986]{ethier1986}
Ethier, S.~N. and Kurtz, T.~G. (1986).
\newblock {\em Markov Processes: Characterization and Convergence}.
\newblock Wiley.

\bibitem[Himmelberg, 1975]{himmelberg1975}
Himmelberg, C.~J. (1975).
\newblock Measurable relations.
\newblock {\em Fund. Math.}, 87:53--72.

\bibitem[Kallenberg, 2002]{kallenberg2002}
Kallenberg, O. (2002).
\newblock {\em Foundations of Modern Probability}.
\newblock Springer, second edition.

\bibitem[Kamae et~al., 1977]{kamae1977}
Kamae, T., Krengel, U., and O'Brien, G.~L. (1977).
\newblock Stochastic inequalities on partially ordered spaces.
\newblock {\em Ann. Probab.}, 5(6):899--912.

\bibitem[Kemeny and Snell, 1960]{kemeny1960}
Kemeny, J. and Snell, J. (1960).
\newblock {\em Finite {Markov} chains}.
\newblock Princeton Universty Press.

\bibitem[Kingman, 1969]{kingman1969}
Kingman, J. F.~C. (1969).
\newblock Markov population processes.
\newblock {\em J. Appl. Probab.}, 6:1--18.

\bibitem[Kuratowski and Ryll-Nardzewski, 1965]{kuratowski1965}
Kuratowski, K. and Ryll-Nardzewski, C. (1965).
\newblock A general theorem on selectors.
\newblock {\em Bull. Acad. Polon. Sci. S\'er. Sci. Math. Astronom. Phys.},
  13:397--403.

\bibitem[Lindvall, 1992]{lindvall1992}
Lindvall, T. (1992).
\newblock {\em Lectures on the Coupling Method}.
\newblock Wiley.

\bibitem[L{\'o}pez et~al., 2000]{lopez2000}
L{\'o}pez, F.~J., Mart{\'i}nez, S., and Sanz, G. (2000).
\newblock Stochastic domination and {Markovian} couplings.
\newblock {\em Adv. Appl. Probab.}, 32:1064--1076.

\bibitem[L{\'o}pez and Sanz, 2002]{lopez2002}
L{\'o}pez, F.~J. and Sanz, G. (2002).
\newblock Markovian couplings staying in arbitrary subsets of the state space.
\newblock {\em J. Appl. Probab.}, 39:197--212.

\bibitem[Marshall and Olkin, 1979]{marshall1979}
Marshall, A.~W. and Olkin, I. (1979).
\newblock {\em Inequalities: Theory of Majorization and Its Applications}.
\newblock Academic Press.

\bibitem[Massey, 1987]{massey1987}
Massey, W.~A. (1987).
\newblock Stochastic orderings for {Markov} processes on partially ordered
  spaces.
\newblock {\em Math. Oper. Res.}, 12(2):350--367.

\bibitem[M{\"u}ller and Stoyan, 2002]{muller2002}
M{\"u}ller, A. and Stoyan, D. (2002).
\newblock {\em Comparison Methods for Stochastic Models and Risks}.
\newblock Wiley.

\bibitem[Serfozo, 1999]{serfozo1999}
Serfozo, R. (1999).
\newblock {\em Introduction to Stochastic Networks}.
\newblock Springer.

\bibitem[Shaked and Shanthikumar, 2007]{shaked2007}
Shaked, M. and Shanthikumar, J.~G. (2007).
\newblock {\em Stochastic orders}.
\newblock Springer Series in Statistics. Springer, New York.

\bibitem[Srivastava, 1998]{srivastava1998}
Srivastava, S.~M. (1998).
\newblock {\em A course on {B}orel sets}.
\newblock Springer.

\bibitem[Strassen, 1965]{strassen1965}
Strassen, V. (1965).
\newblock The existence of probability measures with given marginals.
\newblock {\em Ann. Math. Statist.}, 36(2):423--439.

\bibitem[Szekli, 1995]{szekli1995}
Szekli, R. (1995).
\newblock {\em Stochastic Ordering and Dependence in Applied Probability}.
\newblock Springer.

\bibitem[Thorisson, 2000]{thorisson2000}
Thorisson, H. (2000).
\newblock {\em Coupling, Stationarity, and Regeneration}.
\newblock Springer.

\bibitem[Wagner, 1977]{wagner1977}
Wagner, D.~H. (1977).
\newblock Survey of measurable selection theorems.
\newblock {\em SIAM J. Control Optimization}, 15(5):859--903.

\bibitem[Whitt, 1986]{whitt1986}
Whitt, W. (1986).
\newblock Stochastic comparisons for non-{Markov} processes.
\newblock {\em Math. Oper. Res.}, 11(4):608--618.

\bibitem[Winston, 1977]{winston1977}
Winston, W. (1977).
\newblock Optimality of the shortest line discipline.
\newblock {\em J. Appl. Probab.}, 14(1):181--189.

\end{thebibliography}
\end{document}